%%%%%%%%%%%%%%%%THIS IS AN AMS-TEX DOCUMENT%%%%%%%%%%%%%%%%%%%%%%%%%%%%%%

\magnification=\magstep1
\vsize=22truecm
\input amstex

\documentstyle{amsppt}

\TagsOnRight

\def\cite#1{{\rm [#1]}}

\leftheadtext{E. MAKAI, Jr.}

\rightheadtext{FIVE-NEIGHBOUR PACKINGS OF CENTRALLY SYMMETRIC CONVEX DISCS}

\topmatter

\title
Five-neighbour packings of centrally symmetric convex discs. 
\endtitle

\author
Endre Makai, Jr.$^*$
\endauthor

\address
Endre Makai, Jr.,\newline
\indent Alfr\'ed R\'enyi Mathematical Institute,\newline
\indent Hungarian Academy of Sciences,\newline
\indent H-1364 Budapest, P.O. Box 127,\newline
\indent HUNGARY
\endaddress
\email
makai\@renyi.hu, {\rm{http://www.renyi.hu/\~{}makai}}
\endemail

\thanks
$^*$Research (partially) supported by Hungarian National Foundation for
Scientific Research, grant nos. K68398, T046846.
\endthanks



%%%%%%%%%%%%%%%%%%%%%%%%%%%%%%%%%%%%%%%%%%%%%%%%%%%%%%%%%%%%%%%%%%%%%%%%%%

\abstract
In an old paper of the author the thinnest five-neighbour packing of translates
of a convex
disc (different from a parallelogram) was determined. The minimal density was
$3/7$, and was attained for a certain packing of triangles. In that paper it
was announced that for centrally symmetric convex plates (different from a
parallelogram) the analogous
minimal density was $9/14$, and was attained for a certain
packing of affine regular hexagons, and a very sketchy idea of the proof was
given. In this paper we give details of this proof.
\endabstract

%%%%%%%%%%%%%%%%%%%%%%%%%%%%%%%%%%%%%%%%%%%%%%%%%%%%%%%%%%%%%%%%%%%%%%%%%%

\endtopmatter

%%%%%%%%%%%%%%%%%%%%%%%%%%%%%%%%%%%%%%%%%%%%%%%%%%%%%%%%%%%%%%%%%%%%%%%%%%
%%%%%%%%%%%%%%%%%%%%%%%%%%%%%%%%%%%%%%%%%%%%%%%%%%%%%%%%%%%%%%%%%%%%%%%%%%
%%%%%%%%%%%%%%%%%%%%%%%%%%%%%%%%%%%%%%%%%%%%%%%%%%%%%%%%%%%%%%%%%%%%%%%%%%

\document
\heading
\S1. Introduction
\endheading

For notations cf. \cite{M}.

L. Fejes T\'oth \cite{FTL} (1973)
proved that a $5$-neighbour packing of circles in the plane has
(lower)
density at least $\pi {\sqrt{3}}/7$, which is sharp; and that the density of
a five-neighbour
packing of translates of any convex disc different from a parallelogram is at
least some positive constant, here $5$ being sharp,
and conjectured that the sharp constant is $3/7$; 
and that the density of a $6$-neighbour packing of translates of any convex
disc is at least a $1/2$, here $1/2$ being sharp.

L. Fejes T\'oth-N. Sauer \cite{FTS} (1977)
proved that if for a packing of translates of cubes in
${\Bbb{R}}^n$ the number of at most $k$-th neighbours of any cube is at least
$(k+1)(2k+1)^{n-1}+1$, then the (lower) density is positive, and here 
$(k+1)(2k+1)^{n-1}+1$ is sharp. Moreover, they constructed for $n=2$ and $k=2$
and any integer $i \in [16, 23]$ a packing of translates of a square with at
least $i$ at most second neighbours of low density (conjecturably with minimal
density).

The author \cite{Ma} (1985) proved the conjectured value $3/7$ of \cite{FTL}.

G. Fejes T\'oth-L. Fejes T\'oth \cite{FT-FT} (1991)
proved that if there is a five-neighbour
packing of circles in the plane such
that the infimum of the quotients of any two radii is greater than some number
$h_0$ then the (lower) 
density of the packing is positive, and if this infimum equals
$h_0$ then the packing can have density $0$.

G. Fejes T\'oth \cite{FTG} (1981), H. Sachs \cite{Sa} (1986)
and G. Kert\'esz \cite{K} (1994)
proved that if we have in ${\Bbb{R}}^3$ a $10$-neighbour 
packing of unit balls, then it has a (lower)
density at least some positive constant, and here $10$ is sharp.

K. Bezdek-P. Brass \cite{BB} (2003)
proved that if we have in ${\Bbb{R}}^n$ a $2 \cdot 3^{n-1}$-neighbour
packing of translates of a convex body, then its density is positive, and this
number of neighbours is sharp (for a parallelotope).

It follows from the main Theorem 7 of O. R. Musin \cite{Mu} (2006) 
that if we have in ${\Bbb{R}}^4$ a
$19$-neighbour packing of unit balls, then it has a (lower)
density at least some positive constant, and here $19$ is sharp. (This theorem 
has to be combined with the method of the earlier papers, e.g., \cite{FTL} for
${\Bbb{R}}^2$, and \cite{FTG} for ${\Bbb{R}}^3$.)

%%%%%%%%%%%%%%%%%%%%%%%%%%%%%%%%%%%%%%%%%%%%%%%%%%%%%%%%%%%%%%%%%%%%%%%%%%%
%%%%%%%%%%%%%%%%%%%%%%%%%%%%%%%%%%%%%%%%%%%%%%%%%%%%%%%%%%%%%%%%%%%%%%%%%%%
%%%%%%%%%%%%%%%%%%%%%%%%%%%%%%%%%%%%%%%%%%%%%%%%%%%%%%%%%%%%%%%%%%%%%%%%%%%

\heading
\S2. New results
\endheading

We use the definitions and notations of the paper \cite{M}.

First we state the following

%%%%%%%%%%%%%%%%%%%%%%%%%%%%%%%%%%%%%%%%%%%%%%%%%%%%%%%%%%%%%%%%%%%%%%%%%

\proclaim{Proposition}
Let us have a five-neighbour packing of translates of a planar
convex domain~$D$, which is not a parallelogram.
Then the corresponding polygonal subdivision of ${\Bbb{R}}^2$
(whose edges are the segments
connecting the homologous points of the touching pairs of the translates of the
convex domain $D$) consists of polygons with at most six sides.
\endproclaim

%%%%%%%%%%%%%%%%%%%%%%%%%%%%%%%%%%%%%%%%%%%%%%%%%%%%%%%%%%%%%%%%%%%%%%%%%%

This follows from the existence of a Brass angular measure for any centrally
symmetric convex disc $(D-D)/2$ different from a parallelogram (equivalently:
for $D$ not a parallelogram), cf. the exposition in Swanepoel \cite{Sw}, Ch. 7.

%%%%%%%%%%%%%%%%%%%%%%%%%%%%%%%%%%%%%%%%%%%%%%%%%%%%%%%%%%%%%%%%%%%%%%%%%%%%

We will prove the following

\proclaim{Theorem}
Let us have a packing of translates of a centrosymmetric convex domain~$D$,
corresponding to a subdivision of the plane in convex polygons with at most
six sides.
Let $\lambda$ be the average number of neighbours (in a large circle) and let
$d_0'$ denote the density of the thinnest six-neighbour lattice packing of
translates of~$D$.
Then for the density $d'$ in the large circle we have asymptotically the
following lower bounds:
$$
\align
&\left.
\aligned
\frac{d_0'}{2d_0'(6 - \lambda) + \left(\frac32 \lambda - 8 \right)} \ \
&\text{if }\ 4 \leq \lambda \leq 6,\\
\noalign{\smallskip}
\frac{d_0'}{2d_0'(\lambda - 2) + 6}  \ \ &\text{if }\ 3 \leq \lambda \leq 4,
\endaligned
\right\} d_0' \geq \frac78\,,\\
\noalign{\smallskip}
&\ \, \frac{d_0'}{\frac23 d_0'(6 - \lambda) + \frac13 (\lambda -
3)} \ \ \text{if }\ d_0' \leq \frac78.
\endalign
$$
Equality stands for any $3 \leq \lambda \leq 6$ and any $\frac34 \leq
d_0' \leq 1$ e.g.\ for a centrosymmetric hexagon, one of whose diagonals is
parallel to the corresponding sides.
The ratio of the lengths of these sides and this diagonal is $2d_0' - 1$.
\endproclaim

%%%%%%%%%%%%%%%%%%%%%%%%%%%%%%%%%%%%%%%%%%%%%%%%%%%%%%%%%%%%%%%%%%%%%%%%5

\proclaim{Corollary 1}
Let the conditions of the theorem hold.
Then we have for $d'$ asymptotically the following lower bounds:

\newpage

$$
\align
\frac34 \frac{2}{- \frac23 \lambda + 4} &\text{ if } \ \
4 \frac45 \leq \lambda \leq 6,\\
\noalign{\smallskip}
\frac{2}{8 - \lambda} \quad &\text{ if }\ \ 4 \leq \lambda \leq 4 \frac45,\\
\noalign{\smallskip}
\frac12 \qquad &\text{ if } \ \ 3 \leq \lambda \leq 4.
\endalign
$$
Equality stands for any $4\frac45 \leq \lambda \leq 6$ for a regular hexagon,
for any $3 \leq \lambda \leq 4 \frac45$ for a square.
\endproclaim

%%%%%%%%%%%%%%%%%%%%%%%%%%%%%%%%%%%%%%%%%%%%%%%%%%%%%%%%%%%%%%%%%%%%%%%

\proclaim{Corollary 2}
Let the conditions of the theorem hold, except that $D$ should be an
arbitrary, not necessarily centrosymmetric convex domain.
Let $d$ denote the density of our packing, and $d_0$ the density of the
thinnest six-neighbour lattice packing of translates of~$D$.
Then we have for $\frac{d}{d_0}$ the following lower bounds:
$$
\align
\frac{2}{8 - \lambda}\ \ & \text{ if }\ \ 4 \leq \lambda \leq 6,\\
\frac12\quad \ & \text{ if }\ \ 3 \leq \lambda \leq 4.
\endalign
$$
Equality stands for any $3 \leq \lambda \leq 6$ for a square.
\endproclaim

%%%%%%%%%%%%%%%%%%%%%%%%%%%%%%%%%%%%%%%%%%%%%%%%%%%%%%%%%%%%%%%%%%%%%%%%%%%%

The deduction of the corollaries from the theorem is straightforward, one has
to minimize the corresponding expressions from the theorem.
Further evidently $\frac{d}{d_0} = \frac{d'}{d_0'}$, where $d'$ and $d_0'$ are
the corresponding densities for the centrosymmetrized of~$D$ (i.e., $(D-D)/2$).
The corollaries answer questions of L. Fejes T\'oth.

%%%%%%%%%%%%%%%%%%%%%%%%%%%%%%%%%%%%%%%%%%%%%%%%%%%%%%%%%%%%%%%%%%%%%%%%%%%%
%%%%%%%%%%%%%%%%%%%%%%%%%%%%%%%%%%%%%%%%%%%%%%%%%%%%%%%%%%%%%%%%%%%%%%%%%%%%
%%%%%%%%%%%%%%%%%%%%%%%%%%%%%%%%%%%%%%%%%%%%%%%%%%%%%%%%%%%%%%%%%%%%%%%%%%%%

\heading
\S3. Proofs of the new results
\endheading

{\it Proof\/} of the theorem.
Let first $D$ be a centrosymmetric hexagon of the theorem.
For $\lambda = 6, 4, 3$ we have the following packings.

If $3 \leq \lambda \leq 6$ is arbitrary, we take two of the given three
packings, and on some part of the plane we take one of these packings, and on
the other part the other one.
The ratio of the areas of the parts is chosen so as to assure the given value
of~$\lambda$.
Let $d_0' \geq \frac78$.
If $4 \leq \lambda \leq 6$, we use the first and second packings, while if
$3 \leq \lambda \leq 4$ we use the second and third ones.
If $d_0' \leq \frac78$ we use the first and the third packings.

The parts of the planes, on which the different packings are taken, can be
chosen to consist of parallel strips.
Between any two neighbouring strips we insert rows of 

\newpage

domains, each domain in
the row touching its two neighbours.
Let the density of the inserted rows on the whole plane be~$0$.
If the third packing is considered, in the first and last rows, to each two
domains, touching each other  a third one is attached from right, touching one
of the domains (see the figure).
The parallel strips can be translated with respect to each other in their own
direction arbitrarily, and the rows between them can be translated almost
arbitrarily.
One sees easily that in this way we obtain a packing, to which there
corresponds a subdivision of the plane in convex polygons with numbers of
sides $\leq 6$.

One easily checks that thus we obtain packings satisfying the conditions of
the theorem and having the prescribed densities.
Thus we have equality in the case of the centrosymmetric hexagon of the theorem.

Now we prove for the density $d'$ of a packing satisfying the conditions of
the theorem the required lower bounds.

$D$ is the unit circle of a Minkowski geometry.
Unless the contrary is stated, the lengths of vectors will be measured in this
Minkowski geometry.

We introduce some notations.
$A(D)$ denotes the area of~$D$.
$\Delta$ denotes the maximal area of a triangle with unit sides (in the
Minkowski geometry).
$F_k'(D)$ denotes the maximal area of a convex $k$-gon of unit sides, $k \geq
3$.
$F_k'(\Delta, A(D))$ is the maximum of $F_k'(D)$, when $D$ runs over all
centrosymmetric convex domains, for which $\Delta$ and $A(D)$ have the
prescribed values.
We have $d_0' = \frac12 \frac{A(D)}{4\Delta}$, $\frac34 \leq d_0' \leq 1$.

Evidently $F_3'(\Delta, A(D)) = \Delta$.
The following lemmas are concerned with the values of $F_k'(\Delta, A(D))$ for
$k = 4, 5, 6$.

%%%%%%%%%%%%%%%%%%%%%%%%%%%%%%%%%%%%%%%%%%%%%%%%%%%%%%%%%%%%%%%%%%%%%%%%%%

\proclaim{Lemma 1}
Let a $k$-gon with unit sides be bounded by a simple closed curve.
Then its area is at most $F_k'(D)$.
\endproclaim

%%%%%%%%%%%%%%%%%%%%%%%%%%%%%%%%%%%%%%%%%%%%%%%%%%%%%%%%%%%%%%%%%%%%%%%%

\demo{Proof}
If the $k$-gon, $Q_1 \dots Q_k$, say, is not convex, let us consider a
supporting line of form $Q_i Q_j$, such that one of the open arcs
$\overset{\frown}\to {Q_iQ_j}$ is in the interior of the convex hull of
$Q_1\dots Q_k$.
Let us replace one of the arcs $\overset{\frown}\to {Q_iQ_j}$ of the $k$-gon
by its mirror image through the middle-point of $Q_iQ_j$.
The new $k$-gon is bounded by a simple closed curve, and compared with the
original $k$-gon, it has a greater area and the set of the side-vectors is the
same, only their order is changed.
Thus a finite number of such steps will lead to a convex $k$-gon of unit sides
having a greater area than $Q_1 \dots Q_k$.
Thus our statement follows.
$\blacksquare $
\enddemo

%%%%%%%%%%%%%%%%%%%%%%%%%%%%%%%%%%%%%%%%%%%%%%%%%%%%%%%%%%%%%%%%%%%%%%%%%%%

\proclaim{Lemma 2}
We have $F_5'(\Delta, A(D)) \leq \frac12\bigl[F_4(\Delta, A(D)) + F_6'(\Delta,
A(D))\bigr]$.
\endproclaim

%%%%%%%%%%%%%%%%%%%%%%%%%%%%%%%%%%%%%%%%%%%%%%%%%%%%%%%%%%%%%%%%%%%%%%%%%

\demo{Proof}
We show $F_5'(D) \leq \frac12 \cdot [F_4'(D) + F_6'(D)]$, from which the
statement follows.
Let us have a convex pentagon $P_1\dots P_4$ of unit sides.
The diagonal $P_1 P_3$ cuts the 

\newpage

pentagon in two parts.
By adding to any of these parts its mirror image through the middle-point of
$P_1P_3$, we obtain a quadrangle, and a hexagon of unit sides, both having an
area twice larger than the area of the corresponding part of the pentagon.
Hence by Lemma~1 our statement follows.
$\blacksquare $
\enddemo

%%%%%%%%%%%%%%%%%%%%%%%%%%%%%%%%%%%%%%%%%%%%%%%%%%%%%%%%%%%%%%%%%%%%%%%%%%

\proclaim{Lemma 3}
Let $k$ be even.
Then among the convex $k$-gons with unit sides of maximal area there is a
centrosymmetric one.
\endproclaim

%%%%%%%%%%%%%%%%%%%%%%%%%%%%%%%%%%%%%%%%%%%%%%%%%%%%%%%%%%%%%%%%%%%%%%%%%%

\demo{Proof}
Let us have a convex $k$-gon $P_1\dots P_k$ of unit sides.
The diagonal $P_1 P_{1 + \frac{k}2}$ cuts $P_1\dots P_k$ into two parts.
We consider the part of greater (or equal) area and add to it its mirror image
through the middle-point of $P_1P_{1 + \frac{k}2}$.
Thus we obtain a centrosymmetric $k$-gon of unit sides, $Q_1\dots Q_k$, say,
bounded by a simple closed curve, and having an area not less than $P_1 \dots
P_k$.
By performing all steps of the procedure of Lemma~1, the consecutive pairs of
steps corresponding to supporting lines which are mirror images through the
centre, we obtain a centrosymmetric convex $k$-gon of unit sides, of area not
less than of $P_1\dots P_k$.
$\blacksquare $
\enddemo

%%%%%%%%%%%%%%%%%%%%%%%%%%%%%%%%%%%%%%%%%%%%%%%%%%%%%%%%%%%%%%%%%%%%%%%%%%

\proclaim{Lemma 4}
We have $F_6'(\Delta, A(D)) = A(D)$.
The equality $F_6'(D) = A(D)$ stands iff $D$ is a centrosymmetric hexagon.
\endproclaim

%%%%%%%%%%%%%%%%%%%%%%%%%%%%%%%%%%%%%%%%%%%%%%%%%%%%%%%%%%%%%%%%%%%%%%%%

\demo{Proof}
Let $P_1 \dots P_6$ be a convex hexagon of unit sides.
We may suppose it centrosymmetric by Lemma~3.
We may suppose its three consecutive side-vectors are $(0, -1), (1, 0),
(x,y)$, where $x,y \geq 0$.
Then the area of $P_1 \dots P_6$ will be $1 + x + y$.

On the other hand, $D$ contains the vectors $(0, \pm 1)$, $(\pm 1, 0)$, $(\pm
x, \pm y)$, hence its area is not less than the area of their convex hull,
which is $1 + x + y$.
Hence $F_6'(D) \leq A(D)$, with equality iff $D$ is the convex hull of the
above six vectors, i.e., $D$ is a centrosymmetric hexagon and $P_1$ the
side-vectors of $P_1\dots P_6$ are the radius vectors from the centre to the
vertices of~$D$.
$\blacksquare $
\enddemo

%%%%%%%%%%%%%%%%%%%%%%%%%%%%%%%%%%%%%%%%%%%%%%%%%%%%%%%%%%%%%%%%%%%%%%%%%%%

\proclaim{Lemma 5}
We have $F_4'(\Delta, A(D)) = A(D) - 4\Delta$.
The equality $F_4'(D) = A(D) - 4\Delta$ stands e.g.\ if $D$ is a
centrosymmetric hexagon of the theorem $\left(d_0'
= \frac12 \frac{A(D)}{4\Delta}\right)$.
\endproclaim

%%%%%%%%%%%%%%%%%%%%%%%%%%%%%%%%%%%%%%%%%%%%%%%%%%%%%%%%%%%%%%%%%%%%%%%%%

\demo{Proof}
If $D$ is a centrosymmetric hexagon of the theorem, the quadrangles on the
second figure, constituting a subdivision of the plane, have side-lengths $2$
and area $4(A(D) - 4\Delta)$.

If we estimate the area of a convex quadrangle of unit sides from above, we
may suppose by Lemma~3 that it is a parallelogram.
We show $F_4'(D) \leq A(D) - 4\Delta$, i.e.\ $A(D) \geq F_4'(D) + 4\Delta$.
It suffices to show $A(D) \geq A(P) + 4A(T)$, where $P$ and $T$ are a
parallelogram, and a triangle of unit sides ($A$ denotes area).
$D$ contains the convex hull of the side-vectors of the triangle and of the
parallelogram, and of their negatives, the end-points of all these vectors
being on the  circumference of~$D$.
Let $D_1$ denote the convex hull of these vectors.
We have $D \supset D_1$.

\newpage

We  distinguish two cases.
Let at first between two neighbouring side-vectors of $P$ be no side-vector
of~$T$.
In this case the inequality $A(D) \geq A(D_1) \geq A(P) + 4A(T)$ is evident.

Now let any two neighbouring side-vectors of~$P$ be separated by side-vectors
of~$T$ or by their negatives.
Let $O$ denote the centre of $D$, $T_1, \dots, T_6$ the end-points of the
side-vectors of~$T$, and their negatives, and $P_1, \dots, P_4$ the end-points
of the side-vectors of~$P$, and their negatives (all vectors measured
from~$O$). Let $OP_1$ lie between $OT_1$ and $OT_2$, and $OP_2$ lie between
$OT_3$ and $OT_4$.
Fixing the vectors $OT_i$ we vary the vectors $OP_i$, not decreasing $A(P)$
and not increasing $A(D_1)$.

Let us draw a parallel to $P_2 P_4$ through~$P_1$.
Passing on this parallel in the direction
$\overset\longrightarrow\to{P_4P_2}$, we reach either the segment $T_1T_2$, or
the elongation of the segment $T_3T_2$ beyond~$T_2$.
The point, where we reach one of these lines, is denoted by~$P_1'$.
If $P_1'$ lies on the segment $T_1 T_2$, we pass further on the segment
$P_1'T_2$ until $T_2$ is reached.
In this case let $P_1'' = T_2$, in the other case let $P_1'' = P_1'$.
Similar action is done symmetrically for $P_3'$, and we obtain thus~$P_3''$.
In course of this motion the area of the parallelogram is not decreasing,
while the area of the convex hull of the points $T_i$ and of the vertices of
the parallelogram is not increasing, and all of these points lie on the
boundary of the convex hull.
Let $D_2$ denote the convex hull of the points $T_i$ and $P_1'', P_2, P_3'',
P_4$.
Thus $A(D_1) \geq A(D_2)$ and $A(P_1 P_2 P_3 P_4) \leq A(P_1'' P_2 P_3'' P_4)$.
Starting with the points $T_i$ and $P_1'', P_2, P_3'', P_4$ we pass similarly
from $P_2$ to a new point $P_2''$, and from $P_4$ to $P_4''$.
Let $D_3$ denote the convex hull of the points $T_i$ and $P_1'', P_2'', P_3'',
P_4''$.
Then $A(D_2) \geq A(D_3)$ and $A(P_1'' P_2 P_3'' P_4) \leq A(P_1'' P_2'' P_3''
P_4'')$.

Thus $P_1''$ lies on the elongation of the segment $T_3 T_2$ beyond $T_2$, and
similar statements hold for the other points $P_i''$ too.
We have
evidently $A(D_3) = \frac12 A(P_1'' P_2'' P_3'' P_4'') + 4A(T)$, whence by the
above chains of inequalities
$A(D) \geq A(D_1) \geq A(D_2) \geq A(D_3) = \frac12 A(P_1'' P_2'' P_3'' P_4'')
+ 4A(T) \geq \frac12 A(P_1'' P_2 P_3'' P_4) + 4A(T) \geq \frac12 A(P_1 P_2 P_3
P_4) + 4 A(T) = A(P) + 4A(T)$.
$\blacksquare $
\enddemo

%%%%%%%%%%%%%%%%%%%%%%%%%%%%%%%%%%%%%%%%%%%%%%%%%%%%%%%%%%%%%%%%%%%%%%%%%%

\demo
{\it Proof\/} of the theorem continued.
$\lambda$ is the average number of neighbours, so by Euler's theorem the
average number of the sides of the polygons of the subdivision is
$\sim \frac{\lambda}{\frac{\lambda}{2} - 1}$.
Also by Euler's theorem the quotient of the number of vertices and the number
of polygons of the polygonal subdivision is $\sim \frac1{\frac{\lambda}{2} -
1}$.
Thus we have for the density $d'$ of our packing asymptotically
$$
\align
d' &= \frac{\text{number of vertices}}{\text{number of
polygons}} \cdot \frac{A(D)}{\text{average area of the polygons}} \geq\\
&\geq \frac{1}{\frac{\lambda}{2} - 1} \cdot \frac{A(D)}{\text{\rm conc}\,
F_k'(\Delta, A(D)) \left(\frac{\lambda}{\frac{\lambda}{2} - 1} \right)},
\endalign
$$

\newpage

where $\text{\rm conc}\, F_k'(\Delta, A(D))$ means the concave hull of
$F_k'(\Delta, A(D))$ for $3 \leq k \leq 6$.
\big(We take the value of this function in $\frac{\lambda}{\frac{\lambda}{2} -
1}$.\big)
Here we used that the area of a $k$-gon of the subdivision is $\leq
F_k'(\Delta, A(D))$ and we used Jensen's inequality.

By the remark just before Lemma~1 and in Lemmas 2, 4, 5 we have determined
$\text{\rm conc}\, F_k'(\Delta, A(D))$.
Reminding $d_0' = \frac12 \frac{A(D)}{4\Delta}$, and calculating the last
expression we obtain the formulae for the lower bound of $d'$, given in the
theorem.
We note yet that for $d_0' \leq \frac78$ we have $F_4'(\Delta,
A(D)) \leq \frac13 F_3'(\Delta, A(D)) + \frac23 F_6'(\Delta(A(D))$, thus
$\text{\rm conc}\, F_k'(\Delta, A(D))$
coincides with the concave hull of $F_3'(\Delta, A(D))$, $F_6'(\Delta, A(D))$,
in $k = 3,6$, while for $d_0' \geq \frac78$ we have
$$
F_4'(\Delta, A(D)) \geq \frac13 F_3'(\Delta, A(D)) + \frac23 F_6'(\Delta, A(D)).
$$
$\blacksquare $
\enddemo

%%%%%%%%%%%%%%%%%%%%%%%%%%%%%%%%%%%%%%%%%%%%%%%%%%%%%%%%%%%%%%%%%%%%%%%%%%

\demo{\bf{Remark}}
One might be tempted to think that Corollary~1 and Corollary~2 can be obtained
by determining or estimating from above $\max \frac{F_k'(D)}{A(D)}$, and
$\max \frac{F_k'(D)}{\Delta}$, $3 \leq k \leq 6$.
Evidently we have $\max \frac{F_k'(D)}{A(D)} = \max \frac{F_k'(\Delta,
A(D))}{A(D)}$ and $\max \frac{F_k'(D)}{\Delta} = \max \frac{F_k'(\Delta,
A(D))}{\Delta}$.
Hence $\max \frac{F_k'(D)}{\Delta}$ is attained e.g.\ for a square, for $k =
3,4,6$, and we can obtain in this way $\min \frac{d'}{d_0'}$, as function
of~$\lambda$, as in Corollary~2.
But $\max \frac{F_k'(D)}{A(D)}$ is attained for $k = 3$ iff $D$ is an affine
regular hexagon, and for $k = 4$ iff $D$ is a parallelogram, for $k = 6$ iff
$D$ is a centrosymmetric hexagon, so using this we do not obtain for $d'$ (as
function of $\lambda$) the exact lower bound (except for $3 \leq \lambda \leq
4$) and $\lambda = 6$.
Besides that our theorem is a common generalization of the two corollaries,
this fact lead us to the investigation of $F_k'(\Delta, A(D))$.
\enddemo

%%%%%%%%%%%%%%%%%%%%%%%%%%%%%%%%%%%%%%%%%%%%%%%%%%%%%%%%%%%%%%%%%%%%%%%5

\subheading
{Acknowledgement} The author thanks L. Fejes T\'oth for posing this
problem, and for his encouragement during the work.

%%%%%%%%%%%%%%%%%%%%%%%%%%%%%%%%%%%%%%%%%%%%%%%%%%%%%%%%%%%%%%%%%%%%%

%\subheading
%{References}

\enddocument